\input amstex
\documentstyle{amsppt}
\document
\magnification=1200
\NoBlackBoxes
\nologo
\pageheight{18cm}


\bigskip

\centerline{\bf MODULI STACKS $\overline{L}_{g,S}$}

\medskip

\centerline{\bf Yu.~I.~Manin}

\medskip

\centerline{\it Max--Planck--Institut f\"ur Mathematik, Bonn}

\smallskip

\centerline{\it and Northwestern University, Evanston, USA}

\bigskip

{\bf Abstract.} This paper is a sequel to the paper
by A.~Losev and Yu.~Manin [LoMa1], in which new moduli stacks
$\overline{L}_{g,S}$ of pointed curves were introduced. They classify
curves endowed with a family of smooth points divided into two groups,
such that the points of the second group are allowed to
coincide. The homology of these stacks form components
of the extended modular operad whose combinatorial
models are further studied in [LoMa2].
In this paper the basic geometric properties of 
$\overline{L}_{g,S}$ are established using the notion
of weighted stable pointed curves introduced recently
by B.~Hassett. The main result is a generalization
of Keel's and Kontsevich -- Manin's theorems on the structure
of $H^*(\overline{M}_{0,S}).$

\bigskip

\hfill{{\it \`A Pierre Cartier, en t\'emoignage de respect et d'amiti\'e}}

\medskip

\centerline{\bf \S 0. Introduction}

\medskip

This paper, together with [LoMa2], is a sequel to [LoMa1]
where new moduli stacks $\overline{L}_{g,S}$ were first
introduced. Briefly, let $S$ be a finite set of labels
partitioned into two subsets, white and black ones.
The stack $\overline{L}_{g,S}$ parametrizes 
algebraic curves of genus $g$ endowed with a family
of smooth points labeled by $S$ and satisfying
a certain stability condition. All points endowed with
white labels must be pairwise distinct, and distinct
from all points with black labels. Black points may cluster
in an arbitrary way. If there are no black points,
we get the classical Deligne--Mumford stacks $\overline{M}_{g,S}.$
The simplest non--trivial family of $\overline{L}_{g,S}$
thoroughly studied in [LoMa1] corresponds to the case $g=0$,
two white points, and arbitrary number of black points.
In particular, in [LoMa1] we calculated the  Chow/homology groups
of these stacks, and studied the representation theory
of a family of rings constructed from these groups.

\smallskip

One objective of [LoMa1] was the extension of the homology
operad $\{H_*(\overline{M}_{0,n+1})\}$ by a missing
$n=1$ term. We argued that the union of (the homology of)
all $\overline{L}_{0,S}$ with two white points
provides such a term. The spaces $\overline{L}_{0,S}$
with arbitrary number of white and black points, and more generally,
the stacks $\overline{L}_{g,S}$, are necessary to make the whole
system closed with respect to the clutching
(operadic) morphisms.

\smallskip

In this paper I review the general properties of
$\overline{L}_{g,S}$ and calculate their Chow and (co)homology groups
for $g=0$ using a further generalization of $\overline{M}_{g,S}$
introduced in the recent paper by B.~Hassett [H].
Hassett considers labeling sets $S$ whose elements are endowed with rational
weights, introduces the appropriate notion of stability,
and extends the basic constructions of the theory
of the stacks $\overline{M}_{g,S}$ to the weighted sets.
Our stacks $\overline{L}_{g,S}$  are special cases
of Hassett's stacks corresponding to appropriate weight systems.
Hassett's technics provide a very flexible algebraic
geometric tool and ideally serve our needs.

\smallskip

The paper [LoMa2] is dedicated to the combinatorial models
of the cohomology and homology of $\overline{L}_{0,S}$
and operadic aspects of this
family of rings/modules. The main objective of this paper
is to present a proof that these combinatorial models are actually
isomorphic to the respective cohomology/homology.

\smallskip

The paper is organized as follows. After reminding
principal definitions in \S 1, we introduce the basic structure
morphisms between the stacks $\overline{L}_{g,S}$ in \S 2.
They include {\it clutching morphisms} (translating into
operadic composition laws), {\it forgetting morphisms},
and {\it repainting white to black morphisms}. The latter play
a central role in \S 3 where the main results
of this paper are presented. They concern the structure 
of the Chow/(co)homology groups of $\overline{L}_{0,S}$.

\smallskip

Recall that  the cohomology rings of $\overline{M}_{0,S}$
were calculated by S.~Keel ([Ke]). Keel proved
that $H^*(\overline{M}_{0,S})$ is a quadratic algebra generated by
explicit generators (classes of codimension one
boundary strata indexed by 2--partitions of $S$) 
which satisfy a complete system of explicit linear
and quadratic relations having a transparent geometric origin.
This result was completed in [KM] and [KMK] where a system
of additive generators (classes of boundary strata of arbitrary codimension)
and a complete system of linear relations between them
was deduced from Keel's theorem. The latter description
was crucial for identifying algebras over the homology
operad with formal Frobenius manifolds.

\smallskip

In \S 3 of this paper and in [LoMa2] both theorems 
are generalized to $\overline{L}_{0,S}$. Proofs are rather long
and use algebraic geometric and heavy combinatorial arguments.
The algebraic geometric part takes the original Keel's
theorem as the base of induction and descends to the more general
$S$ by repainting white points to black one by one. The combinatorial
part consists in the thorough study of the abstract quadratic
algebras generated by the (analogs of) Keel's generators and relations.
One intermediate combinatorial result (Theorem 3.6.1) is taken for granted here;
its proof is contained to [LoMa2].

\bigskip

\centerline{\bf \S 1. Definitions and notation}

\medskip

{\bf 1.1. Curves, pointed curves, stability.} {\it A semi--stable curve $C$} over an algrebraically
closed field $k$ is a proper reduced one--dimensional algebraical
scheme over this field having only ordinary double points as singularities.
(Geometric) genus of $C$ is $g:=\roman{dim}\,H^1(C,\Cal{O}_C).$

\smallskip

Let $S$ be a finite set. {\it An $S$--pointed curve $C$}
is a system $(C, x_s\,|\,s\in S)$ where $x_s$ is a family of closed
non--singular $k$--points of $S$, non necessarily pairwise distinct.
The element $s$ is called the label of $x_s$.

\smallskip

The normalization $\widetilde{C}$ of $C$ is a disjoint union
of smooth proper curves. Each irreducible component of
$\widetilde{C}$ carries inverse images of some labeled points $x_s$
and of singular points of $C$. Taken together, these points are called
{\it special ones.} Instead of passing to the normalization,
we may consider branches (local irreducible germs) of $C$ passing through labeled or
singular points. They are in a natural bijection
with special points. 

\smallskip

{\it A painted finite set $S$} is $S$ together with its partition into two
disjoint subsets $S=W\cup B.$ Labels from $W$ (resp. $B$) are called {\it white}
(resp. {\it black}) ones. We may refer to $x_s$ as a white (resp. black) point, 
if its label $s$ is white (resp. black).
 
\smallskip

In the  remainder of this paper, 
{\it  (sets labeling) distinguished points on curves are
usually assumed to be painted in this sense.} This refers
to the labeled points $x_s$ as well as the branches at singular
points. The latter are always painted white.

The following is the main definition of this section.

\medskip

{\bf 1.1.1. Definition.} {\it A semistable $S$--pointed curve is called 
painted stable,
if it is connected and the following conditions
are satisfied:

\smallskip

(i) Each genus zero component of the normalization $\widetilde{C}$
contains at least $3$ special points of which at least 2 are white.
Each genus one component of $\widetilde{C}$
contains at least one special point.

\smallskip

(ii) $x_s$ satisfy the following
condition:
$$
x_s\ne x_t\quad \roman{for\ each}\ s\ne t,\ s\in W,\ t\in S.
\eqno(1.1)
$$
}

\medskip

Notice that in [LoMa1], 4.3.1, such curves were called stringy stable.

\smallskip

In the intermediate constructions we will also widely use
{\it weights} on $S$, in the sense of Hassett ([H]),
and the related notions of weighted stability.
In fact, our paintings corresponds to special systems of Hassett's
weights, and his general techniques
provide an extremely efficient way of working with
moduli stacks of painted curves.

\medskip

{\bf 1.2. Weights and weighted stability.} Let $S$ be an (unpainted) finite set
of labels.  {\it The weight data} on $S$ is
a function $\Cal{A}:\,S\to \bold{Q}$, $s\mapsto a_s$, $0<a_s\le 1.$
We call $S$ together with a weight data {\it a weighted set.}

\medskip

{\bf 1.2.1. Definition ([H]).} {\it A semistable $S$--pointed curve 
$(C, x_s\,|\,s\in S)$ is called 
weighted stable (with respect to $\Cal{A}$)
if the following conditions
are satisfied:

\smallskip

(i) $\omega_{C}(\sum_s a_sx_s)$ is ample where $\omega_{C}$ is the 
dualizing sheaf of $C$.

\smallskip

(ii) For any subset $S'\subset S$ such that $x_t$ pairwise coincide
for $t\in S'$, we have $\sum_{t\in T} a_t \le 1.$}

\smallskip

In (i) and below, we use a shorthand notation. If at least one $a_s$
is $\ne 1$, then only $\omega_{C}^d(\sum_s da_sx_s)$ is an actual
invertible sheaf where $d>1$ is a common denominator
of all $a_s$. Ampleness refers to any of these sheaves.

\smallskip

Clearly, (i) implies that $2g-2+\sum_s a_s >0.$ 

\medskip

{\bf 1.2.2. Lemma.} {\it Let $S=W\cup B$ be a painted set.
Consider a weight data $\Cal{A}$ on $S$ satisfying the following conditions:

\medskip

(*) $a_s=1$ for all $s\in W$, and $\sum_{t\in B} a_t \le 1.$

\medskip

Then a semistable $S$--pointed curve 
$(C, x_s\,|\,s\in S)$ is painted stable if and only if it is weighted stable
with respect to $\Cal{A}$.}

\smallskip

{\bf Proof.} The lift of $\omega_{C}(\sum_s a_sx_s)$ to any component $D$
of the normalization of $C$ embeds into
$\Omega_{D}^1(\sum_t b_tx_t)$ where summation is taken over
all special points of $D$ and the weight $\Cal{A}$
is extended to singularities by 1. Sections which come from
$C$ are singled out by the local condition: 
at pairs of points which get identified in $C$ the sum of residues
vanishes. This shows that the stability condition (i)
of the definition 1.2.1 is non--empty only on
components of genus 0 and 1. When (*) is satisfied,
this stability condition  is then equivalent to the
condition (i) of the Definition 1.1.1.

\smallskip

Similarly, (*) and 1.2.1(ii) taken together say that on a weighted 
stable curve black points may 
pairwise coincide, but
each white point must be different from all other points.
This is precisely the condition (ii) of painted stability.

\medskip

{\bf 1.3. Families of pointed curves.} Let $T$ be a scheme, $S$ a finite set, $g\ge 0$. 
{\it An $S$--pointed  curve (or family of curves) of genus $g$
over $T$} consists of the data
$$
(\pi:\,C\to T;\,x_i:\,T\to C,\ i\in S)
$$
where

\smallskip

(i) $\pi$ is a flat proper morphism whose geometric fibres $C_t$
are semistable curves of genus $g$.

\smallskip

(ii) $x_i, i\in S,$ are sections of $\pi$ not containing
singular points of geometric fibres.

\smallskip

Various definitions of stability from the subsections 1.2--1.3 are 
generalized to families by
requiring the respective properties to hold on all geometric
fibers of $\pi.$ Ampleness condition can be equivalently stated
in terms of the relative dualizing sheaf.

\medskip

{\bf 1.3.1. Stacks of weighted stable curves $\overline{\Cal{M}}_{g,\Cal{A}}.$} 
The first main result of [H]
is a proof of the following fact. Fix a weighted set of labels $S$
and a value of genus $g$. Then families of weighted stable
$S$--pointed curves of genus $g$ form (schematic points of) a connected smooth proper over $\bold{Z}$
Deligne--Mumford stack $\overline{\Cal{M}}_{g,\Cal{A}}.$ The respective
coarse moduli scheme is projective over $\bold{Z}$.

\medskip

{\bf 1.3.2. Stacks of painted stable curves $\overline{L}_{g,S}.$} Let now 
$S$ be a painted set. Reinterpreting the painted stability
condition as in Lemma 1.2.2, we deduce from Hassett's theorem
the existence of the respective stacks $\overline{L}_{g,S}$
which were constructed in [LoMa1] by an alternative method which we
called there ``adjunction of the generic black point.''

\medskip

{\bf 1.4. Graphs.} {\it A graph $\tau$}, by definition,
is a quadruple $(V_{\tau}, F_{\tau}, \partial_{\tau},j_{\tau})$
where $V_{\tau}$, resp. $F_{\tau}$, are finite sets of
vertices, resp. flags; $\partial_{\tau}:\, F_{\tau}\to V_{\tau}$
is the boundary, or incidence, map; $j_{\tau}:\,F_{\tau}\to F_{\tau}$ is an involution of the set of flags. {\it The geometric realization
of $\tau$}  is a topological space which is obtained from
$F_{\tau}$ copies of $[0,1]$ by gluing together
points $0$ in the copies corresponding to each vertex $v\in V_{\tau}$,
and by gluing together points $1$ in each orbit of $j_{\tau}$.
This motivates considering the following auxiliary sets
and their geometric realizations: the set $E_{\tau}$ of {\it edges} of
$\tau$, formally consisting of cardinality two orbits of $j_{\tau}$,
and the set $T_{\tau}$ of {\it tails}, consisting of those
flags, which are $j_{\tau}$--invariant.

\smallskip

We will mostly think and speak about graphs directly in terms
of their geometric realizations. In particular, $\tau$ will be called
{\it connected} (resp. {\it tree}, resp. {\it forest})
if its geometric realization is connected (resp. connected and
has no loops, resp. is 
a disjoint union of trees).

\smallskip

The (dual) modular graph of an $S$--pointed semistable curve is defined
in the same way as in the usual case. We use the conventions
of [Ma], III.2 where the reader can find further details.
Briefly, irreducible components become vertices,
pairs of special points that are identified become edges,
labeled points become tails, so that
tails acquire labeling by $S$. Tails now can be of two
types, we may refer to them and their marks 
as ``black'' and ``white'' ones as well, and call the graph
painted one. Moreover, each vertex is marked by
the genus of the corresponding (normalized) component. 

\smallskip

In the genus zero case, all relevant graphs are trees.
If we delete a vertex $v$ from the geometric realization of the tree 
$\tau$, it will break into a set of $\ge 3$
connected components which we will call
{\it branches} of $\tau$ at $v$. Their set
is canonically bijective to the set
of flags $F_{\tau}(v)$ incident to $v$: we can say that
the  branch starts with the respective flag. 
In the extreme case, a branch can be a single tail.

\bigskip

\centerline{\bf \S 2. Basic morphisms}

\medskip

{\bf 2.1. Clutching morphisms and boundary strata.} First, recall
a general construction studied in [Kn], \S 3. Let
$C/T$ be a flat family of semistable curves over a scheme $T$
endowed with two non--intersecting sections $x_s,x_t$.
Then there is another family of curves $C'/T$
and a morphism $p:\,C\to C'$ over $T$ such that
$p\circ x_s= p\circ x_t$ and $p$ is universal with this poperty.
Knudsen proves that $p$ is a finite morphism, and $C'/T$
is again semistable. 
\smallskip

We will apply this construction and its iterations
to the (schematic points of) various products
of $\overline{L}_{g,S}$.  In all instances,
we will be gluing together only {\it pairs of white points.}

\smallskip

First, let $S_1=S'\cup \{s\}$, $S_2=S''\cup \{t\}$, 
with $s, t$ white. Put $S:=S'\cup S''$. Then $\overline{L}_{g_1,S_1}
\times \overline{L}_{g_2,S_2}$ carries a family of disconnected
curves $pr_1^*(C_{g_1,S_1})\coprod_T pr_2^*(C_{g_2,S_2})$
endowed with two nonintersecting sections coming from $x_s, x_t$.
Clutching them together, we obtain
a morphism of stacks
$$
\overline{L}_{g_1,S_1}
\times \overline{L}_{g_2,S_2} \to \overline{L}_{g_1+g_2,S}
\eqno(2.1)
$$
which is a closed immersion and defines a boundary divisor
of $\overline{L}_{g,S}$, $g=g_1+g_2$. One can similarly
define boundary divisors obtained by gluing together
a pair of white sections $s,t\in S'$ of the universal
curve over $\overline{L}_{g',S'}$. Such divisors exist only when $g>0$.

\smallskip

More generally, let $\tau$ be a graph whose tails are bijectively
labeled by a painted set $S$ and each vertex $v$
is endowed with a value of ``genus'' $g_v$. All these labels
form a part of the structure of $\tau$. 
We extend the painting to all flags: halves of edges are white.
Assume that  
$\tau$ is {\it painted stable}: vertices of genus $0$ are incident
to $\ge 3$ flags of which at least $2$ are white, and 
vertices of genus 1 are incident to $\ge 1$ flags.
Then we put 
$$
\overline{L}_{\tau}:=\prod_{v\in V_{\tau}} \overline{L}_{g_v,F_\tau (v)}
\eqno(2.2)
$$
where $F_{\tau}(v)$ is the set of flags incident to $v$.
This stack carries the disjoint union of universal
curves lifted from the factors. This union is endowed
with a family of pairwise disjoint sections
corresponding to the halves of all edges of $\tau$.
We can now clutch together pairs of sections corresponding
to the halves of one and the same edge. This produces
a boundary stratum morphism
$$
\overline{L}_{\tau} \to \overline{L}_{g, S}
\eqno(2.3)
$$
where $g=\sum g_v + \roman{rk}\,H_1(\tau ).$ In the genus zero case,
we should consider only trees whose all vertices
have genus zero so that one can forget about the latter labels.

\medskip

{\bf 2.2. Forgetting morphisms.} Let $S$ be a painted set, $g\ge 0$,
and $S'\subset S$ a subset such that $\overline{L}_{g,S}$
and $\overline{L}_{g,S'}$ are nonempty. Then
we have a canonical forgetting morphism
$$
\varphi :\, \overline{L}_{g,S} \to \overline{L}_{g,S'}
\eqno(2.4)
$$
which on the level of curves consists in forgetting
the sections labeled by $S-S'$ and consecutively
contracting the components that become unstable.

\smallskip

This is a particular case of forgetting morphisms
defined by Hassett ([H], Theorem 4.3). Hassett's
theorem is applicable thanks to Lemma 1.2.2.

\medskip

{\bf 2.3. Repainting morphisms.} Let $S$ be a painted set,
$a\in S$ a white label. Denote by $S'$ the painted set
with the same elements in which now $a$ is painted black
whereas all other labels keep their initial colors.
Assuming again that $\overline{L}_{g,S}$
and $\overline{L}_{g,S'}$ are nonempty we have a repainting
morphism
$$
\rho :\, \overline{L}_{g,S} \to \overline{L}_{g,S'}
\eqno(2.5)
$$
which on the level of curves consists in repainting black
the section $x_a$ and again consecutively
contracting the components that become unstable.

\smallskip

This is a particular case of reduction morphisms
defined by Hassett [H], Theorem 4.1 which is
applicable thanks to Lemma 1.2.2 as well.

\bigskip

\centerline{\bf \S 3. Chow groups of $\overline{L}_{0,S}$} 

\medskip

{\bf 3.1. Keel's relations for $A^*(\overline{L}_{0,S})$.} In the following
we will be considering only $S$--ponted painted stable curves of genus zero.
Painted stability of various sets, partitions, trees etc.
means that the respective stacks are nonempty.
For any painted stable 2--partition $\sigma$ of $S$ we denote
by $[D(\sigma )]\in A^1(\overline{L}_{0,S})$ the class of the respective
boundary divisor (2.1). Call an ordered quadruple of pairwise distinct
elements $i,j,k,l\in S$ {\it allowed}, if
both partitions $ij|kl$ and $kj|il$ are painted stable.
For a painted stable $\sigma$ put $\epsilon (\sigma ;i,j,k,l)=1$
if $\{i,j,k,l\}$ is allowed and $ij\sigma kl$; $-1$,
if $\{i,j,k,l\}$ is allowed and $kj\sigma il$;
and $0$ otherwise. The following theorem generalizing
Keel's presentation is the main
result of this section.

\medskip

{\bf 3.1.1. Theorem.} {\it The classes $[D(\sigma )]$ generate the ring
$A^*(\overline{L}_{0,S})$. They satisfy the following relations (3.1), (3.2)
which provide a presentation of this ring. First, for each
allowed quadruple $i,j,k,l$,
$$
\sum_{\sigma} \epsilon (\sigma ;i,j,k,l)\,[D(\sigma )] =0\, .
\eqno(3.1)
$$
Second, let $\sigma,\,\sigma'$ be two stable painted partitions
such that there exists
an allowed quadruple $i,j,k,l$ with $ij\sigma kl,\,kj\sigma' il$.
Then
$$
[D(\sigma )]\,[D(\sigma' )]=0.
\eqno(3.2)
$$
}

\smallskip
We break the algebraic--geometric arguments in the 
proof into a series of Lemmas. 

\medskip

{\bf 3.2. Lemma.} {\it The classes $[D(\sigma )]$ satisfy
(3.1) and (3.2).}

\smallskip

{\bf Proof.} In fact, let $\{i,j,k,l\}$ be allowed.
Consider the forgetting morphism $\varphi:\,
\overline{L}_{0,S}\to \overline{L}_{0,\{ijkl\}}$.
Partitions $ij|kl$ and $kj|il$ define two boundary
points in $\overline{L}_{0,\{ijkl\}}\cong \bold{P}^1.$
Their inverse images are precisely sums of boundary divisors 
of $\overline{L}_{0,S}$
entering in (3.1) with coefficients 1, resp. $-1$.
This gives (3.1) whereas (3.2) follows from the fact
that fibers of $\varphi$ over different points do not intersect.

\medskip

{\bf 3.3. Lemma.} {\it The classes $[D(\sigma )]$ additively
generate $A^1(\overline{L}_{0,S})$.}

\smallskip

{\bf Proof.} This was proved by Keel ([Ke]) for the case when
all labels are white. Consider a painted set
$S$ with $\ge 3$ white labels. Choose a white label
$a\in S$ and repaint it black. Denote the resulting
painted set $S'$. Consider the
repainting morphism $\rho :\,\overline{L}_{0,S}\to
\overline{L}_{0,S'}.$ Since $\rho$ is birational,
$\rho_*:\, A^k(\overline{L}_{0,S})\to A^k(\overline{L}_{0,S'})$
is surjective for each $k$. We will show that $\rho_*([D(\sigma )])$
is a linear combination of boundary divisors for
each painted stable 2--partition $\sigma$ of $S$. This will prove
our statement by induction on the number of black labels.

\smallskip

If $\sigma$ remains painted stable after repainting $a$,
we have simply $\rho_*([D(\sigma )])= [D(\sigma')]$
where $\sigma'$ is the same partition of $S'$.

\smallskip

If $\sigma$ becomes unstable, one part of it
must be $\{a\}\cup F$ where $F\subset S$ consists only of black labels.
When $|F|\ge 2,$ we have $\rho_*([D(\sigma )])=0.$
In fact, according to [H], Prop. 4.5, $[D(\sigma )]$
is contracted by $\rho$.

\smallskip

Finally, assume that one part of $\sigma$ is of the form
$\{a,b\}$ where $b$ is black. Choose two white labels $i,j$ in the other part of 
$\sigma$. The quadruple $i,j,a,b$ is allowed in $S$.
Write the relation (3.1) for it as an expression
for $[D(\sigma )]$:
$$
[D(\sigma )]=-\sum_{\tau :\,ij\tau ab} [D(\tau)] +\sum_{\tau :\,aj\tau ib} [D(\tau)]
\eqno(3.3)
$$
where in the first sum the part containing $a,b$ must contain
at least one more label. It is clear now that
applying $\rho_*$ to any summand in the right hand side
we get either a boundary divisor, or zero. This completes the proof.

\medskip

We will now generalize these results to the classes of boundary strata
of arbitrary codimension. Consider a painted stable $S$--tree
$\tau$, a vertex $v$ of it and an allowed quadruple
of flags $I,J,K,L$ at $v$ (recall that all halves of edges are
painted white). For any painted stable 2--partition $\alpha$ of $I,J,K,L$
define the respective tree $\tau (\alpha )$ which has one extra edge
replacing $v$ in $\tau$ and which breaks $F_{\tau} (v)$ according to $\alpha$
(cf. [LoMa2], 1.3--1.4).

\medskip

{\bf 3.4. Lemma.} {\it For each $(\tau , v;I,J,K,L)$ as above, we have the following
relation between boundary strata in $A^*(\overline{L}_{0,S})$:
$$
\sum_{\alpha} \epsilon (\alpha ;I,J,K,L)\, [D(\tau (\alpha ))] = 0\, .
\eqno(3.4)
$$
}

\smallskip

{\bf Proof.} Notice that when $\tau$ is the one--vertex $S$--tree, (3.4)
reduces to (3.1). Conversely, (3.4) can be deduced from (3.1)
in the following way. The closed stratum $D(\tau )$
is (the image of) $\prod_{w\in V_{\tau}}\overline{L}_{0,F_{\tau}(w)}$
in $\overline{L}_{0,S}$. Replacing in (3.1) the label set $S$ by $F_{\tau}(v)$ and 
$i,j,k,l$ by $I,J,K,L$, we get a relation in $A^*(\overline{L}_{0,F_{\tau}(v)})$.
Tensor multiplying
this identity by the fundamental classes of all remaining
$\overline{L}_{0,F_{\tau}(w)}$ and taking the direct image in $A^*(\overline{L}_{0,S})$
we finally obtain (3.4).

\medskip

{\bf 3.5. Lemma.} {\it  The classes $[D(\tau )]$ for painted stable
$S$--trees $\tau$ additively
generate $A^*(\overline{L}_{0,S})$.}

\smallskip

{\bf Proof.} We extend the proof of Lemma 3.3 to this case,
by starting with Keel's result for the case when all labels are 
white, and repainting the necessary amount of white points one by one.

\smallskip

Again, we have to check only that $\rho_*$ maps classes of boundary strata
to linear combinations of such classes. We may and will assume that
$\tau$ has at least two edges. As above, if $\tau$ remains
stable after repainting $a$, this is clear. If $\tau$ becomes unstable,
the tail $a$ must be incident to an end vertex $v$ of $\tau$,
and all other tails at this vertex must be black. Let their set
be $F$. We will first check
that if the number of these black tails is $|F|\ge 2$,
then $\rho_*([D(\tau )])=0.$ In fact, Hassett's description
of the repainting map shows that 
$\rho (D(\tau ))$  parametrizes curves
in which the structure points $x_s, s\in F,$ all have to coincide
after the component on which they formerly
freely moved in $C_{0,S}$ has been collapsed. Hence 
$\roman{dim}\,\rho (D(\tau ))< \roman{dim}\,D(\tau )$ so that
$\rho_* [D(\tau )]=0.$

\smallskip

Consider now the case when an end vertex $v$ of $\tau$
carries only two flags $a,b$, and $b$ is black.

\smallskip

Thus $v$ carries three flags of which two
are white and one black. Let us generally call 
such a vertex {\it critical one}. There is a unique maximal
sequence of vertices $v:=v_0,v_1,\dots ,v_{n-1},v_n$,
$n\ge 1$, in $\tau$,
with the following properties:

\smallskip

(a). {\it For each $0\le i\le n-1$, $v_i$ and $v_{i+1}$ are 
opposite vertices of an edge $e_i$.}

\smallskip

(b). {\it $v_0,\dots ,v_{n-1}$ are critical vertices
whereas $v_n$ is not critical.}

\smallskip

For $n\ge 2$, if we delete the vertex $v_n$ from $\tau$,
the connected component containing $a$ will be called 
{\it the critical branch of $a$}, and $n-1$ will be called its
length. 
\smallskip

The vertex $w:=v_n$ completing the critical
branch can fail to be critical in one of two ways:

\smallskip

(b1). {\it $F_{\tau}(w)$ contains only two white flags
but $\ge 2$ black flags.}

\smallskip

(b2). {\it $F_{\tau}(w)$ contains at least three white flags.}

\smallskip

I contend that in the case (b1)
we again have $\rho_*[D(\tau )]=0$.  In fact, consider
a curve $C$ corresponding to a generic geometric point
of $D(\tau ).$ Using Hassett's description of the repainting map,
we see that after repainting $a$ black, all components
of $C$ corresponding to $v_0,\dots ,v_n$ get collapsed.
Collapsing critical components $v_0,\dots ,v_{n-1}$
does not diminish the number of moduli since
they are all projective lines with three special points.
However, collapsing the first non--critical component $v_n$
does diminish the number of moduli since it carries
$\ge 4$ speial points.

\smallskip

It remains to treat the case (b2). 

\smallskip

Since $v_n$ carries at least three white flags,
we can choose at $v_n$ two white flags $K,L$
which do not coincide with the half of the edge $e_{n-1}$.
Since $v_{n-1}$ is critical, besides the other half
of $e_{n-1}$ it carries a black tail which we call $I$
and one more white flag which we call $j$.
Denote by $\sigma$ the result of collapsing $e_{n-1}$ into a vertex $u$
and write the relation (3.4) for $(\sigma ,u;I,J,K,L):$
$$
\sum_{\alpha:\,IJ\alpha KL} [D(\sigma (\alpha ))] -
\sum_{\alpha:\,KJ\alpha IL} [D(\sigma (\alpha ))] =0.
\eqno(3.5)
$$
We will keep denoting $e_{n-1}$ the edge replacing $u$
in $\sigma (\alpha )$, and $v_{n-1}$,
$v_n$ its respective vertices.
In the first sum, there is one term $[D(\tau )]$
corresponding to the partition $\alpha =IJ|...$
of $F_{\sigma}(u)$.  For all other terms,
$v_{n-1}$ will cease to be a critical vertex
since it will carry $\ge 4$ flags.
In the second sum, $v_{n-1}$ is never critical, because
it carries a white flag $K$ and two white halves
of the edges $e_{n-1}$ and $e_{n-2}$
(if $n=1$, in place of $e_{n-2}$ we have $a$).

\smallskip

Thus, (3.5) allows us to replace the Chow class
of $D(\tau )$ by a linear combination of classes
whose critical branches (of $a$) are shorter than that of $\tau$.
Applying this procedure to each term of this expression
we can reduce the length  once more if need be.
But when the length becomes zero, the vertex
$v_0$ ceases to be critical, and the repainting
of the respective class was described at the beginning. 
This completes the proof of Lemma 3.5.

\smallskip

The remaining part of the proof of Theorem 3.1.1 is purely combinatorial.
It relies upon a theorem on the structure of the abstract ring
$H^*_S$ whose presentation is given by (3.1) and (3.2). This theorem
is proved in [LoMa2]. Below we will summarize the necessary information
from [LoMa2] and then complete the proof.

\medskip

{\bf 3.6. The ring  $H^*_S$.} Let
$S$ be a painted set with $|S|\ge 3$ containing at least two
white elements, $k$ a commutative coefficient ring.
Consider the family of independent commuting variables
$\{l_{\sigma}\}$ indexed by painted stable unordered
2--partitions $\sigma$ of $S$ and put $\Cal{R}_S:=k[l_{\sigma}].$

\smallskip

For an allowed quadruple $i,j,k,l\in S$, put
$$
R_{ijkl}:=\sum_{\sigma} \epsilon (\sigma ;i,j,k,l)\,l_{\sigma} \in \Cal{R}_S ,
\eqno(3.6)
$$
For two partitions $\sigma ,\sigma'$ such that there exists
an allowed quadruple $\{i,j,k,l\}$ with $ij\sigma kl,\,kj\sigma il$, put
$$
R_{\sigma\sigma'}:=l_{\sigma}l_{\sigma'}.
\eqno(3.7)
$$
Denote by $I_S\subset \Cal{R}_S$ the ideal generated by all
elements (3.6) and (3.7) and define {\it the combinatorial cohomology
ring} by
$$
H_S^*:=\Cal{R}_S/I_S.
$$
\smallskip

If $\tau$ is a painted stable $S$--tree, we put
$$
m(\tau ):=\prod_{e\in E_{\tau}}l_{\sigma_e} \in \Cal{R}_S
$$
where $\sigma_e$ for an edge $e$ of $\tau$ denotes the
2--partition of $S$ obtained by cutting $\tau$ in a midpoint of $e$.
Monomials $m(\tau )$ are called {\it good}. Such a monomial depends
only on the $S$--isomorphism class of $\tau$. If $\tau$ is one--vertex
tree, we put $m(\tau )=1.$ For any $m\in \Cal{R}_S$, we put $[m]:=m\,\roman{mod}\,I_S\in H^*_S.$

\smallskip

In [LoMa2], the following result is proved:

\medskip

{\bf 3.6.1. Theorem.} {\it (i)  $H^*_S$ as a $k$--module is spanned by the classes
of good monomials $[m(\tau )].$

\smallskip

(ii) Let $(\tau ,v;I,J,K,L)$ run over systems described before Lemma 3.4.
Each such  system determines  the following relation
between good monomials 
(notation being as in (3.4)):
$$
R(\tau ,v; I,J,K,L):=\sum_{\alpha} \epsilon (\alpha ;I,J,K,L)\, m(\tau (\alpha )) \in I_S\, .
\eqno(3.8)
$$
Moreover, (3.8) span
all relations between the classes of good monomials in $H^*_S$.
}
\medskip

Combining this result with Lemmas 3.2--3.5, we get the following statement:

\medskip

{\bf 3.6.2. Proposition.} {\it The map $[l_{\sigma}]\mapsto [D(\sigma )]$ extends
to a surjective ring homomorphism $h_S:\,H^*_S\to A^*(\overline{L}_{0,S})$,
which sends $[m(\tau )]$ to $[D(\tau )]$ for each painted stable
$S$--tree $\tau$.}

\medskip

{\bf 3.6.3. Relations between relations.} For further reference, I collect here
several relations between the elements (3.8):
$$
R(\tau ,v; I,J,K,L) =  R(\tau ,v; K,L,I,J)= 
$$
$$
=R(\tau ,v; J,I,L,K)= -R(\tau ,v; I,L,K,J)\, ,
$$
$$
R(\tau ,v; I,J,K,L) =R(\tau ,v; I,J,L,K) +R(\tau ,v; I,K,J,L) \, ,
\eqno(3.9)
$$
$$
R(\tau ,v; I,J,K,L) =R(\tau ,v; M,K,J,I)+R(\tau ,v; M,I,L,K) \, .
\eqno(3.10)
$$
The first group of relations is straightforward. In (3.10), we assume that
$M\notin \{I,J,K,L\}$ is an extra flag in $F_{\tau}(v)$ such that all involved
quadruples are allowed. To check (3.10), it is convenient to use the following shorthand
notation: denote e.~g. by $LMK|JI$ the sum of those terms in (3.8)
for which $L,M,K$ get into one part of $\alpha$, whereas $J,I$ get into another part. 
Hence we have, for example, $IJ|KL =IJM|KL + IJ|KLM$. Then
$$
R(\tau ,v; M,K,J,I) = MK|JI - MI|JK=
$$
$$
=LMK|JI +MK|JIL -LMI|JK - MI|JKL\, ,
$$
$$
R(\tau ,v; M,I,L,K) = MI|LK - MK|LI = 
$$
$$
=MI|LKJ + MI|LKM -MKJ|LI -MK|LIJ\, .
$$
After adding up and canceling, stop tracking where $M$ goes. We get $IJ|KL - KJ|IL$
which is $R(\tau ,v; I,J,K,L)$.

\medskip

{\it Questions.} Are all additive relations between relations generated by 
(3.9) and (3.10)? Do (3.9) -- (3.10) form the beginning of an interesting
resolution of $H^*_S$? Are  $H^*_S$ Koszul quadratic algebras?

\medskip

{\bf 3.7. End of the proof of Theorem 3.1.1: the strategy.} It remains to establish that
$h_S$ is injective. This was proved by Keel in the case when
all labels are white. We will again  argue by induction on the
number of black points. Consider a painted set $S$
with $\ge 3$ white labels and fix once for all a white label $a\in S$.  Let $S'$ be obtained from $S$
by repainting $a$ to black one.  Using a long and
convoluted inductive procedure we will construct a surjective graded homomomorphism of 
$k$--modules $\rho^H_*:\,H^*_S \to H^*_{S'}$
such that the following diagram  commutes:
$$
\CD
 H_S^* @>h_S>> A_{S}^*\\
@V\rho^H_* VV   @VV\rho_*V \\
H_{S'}^* @>h_{S'}>> A_{S'}^*\\
\endCD
\eqno(3.11)
$$
where $A^*_S:=A^*(\overline{L}_{0,S})$ with coefficients in $k$
(or the cohomology, which is the same), and $\rho_*$ is induced
by the repainting morphism.
By induction, we can assume that $h_S$ is an isomorphism.
Hence to show that $h_{S'}$ is an isomorphism it suffices to check that
$\roman{Ker}\,h_{S'}=0$. To this end we will use
the inverse image repainting morphism
$\rho_H^*:\,H_{S'}^*\to H_S^*$ fitting into another
commutative diagram
$$
\CD
 H_{S'}^* @>h_{S'}>> A_{S'}^*\\
@V\rho_H^* VV   @VV\rho^*V \\
H_S^* @>h_S>> A_S^*\\
\endCD
\eqno(3.12)
$$
and such that $\rho^H_*\circ \rho_H^* =\roman{id}.$ Its construction
is quite easy: see 3.9.1 below.

\smallskip 

When this is achieved, the completion of the proof is
straightforward. Namely, let $\eta \in H_{S'}^*$ be such that
$h_{S'}(\eta )=0.$ Since $\rho^*$ is injective and $h_S$
is an isomorphism, from (3.12) we get $\rho^*_H(\eta )=0$,
and then $\eta =\rho_*^H\circ\rho^*_H(\eta )=0.$

\medskip

{\bf 3.8. Construction of $\rho^H_*$: the induction parameter
and the inductive statements.} 
Our definition of $\rho_*^H(m(\tau))$ is motivated by the calculation
of $\rho_*$ in the proofs of 3.3 and 3.5. We will do it
by induction on the value of the function
$$
l(\tau )=l(\tau ,a) : =\ \roman{the\ length\ of\
the\ critical\ branch\ of}\ \tau \, .
$$
Recall that a critical vertex is a vertex carrying exactly
two white flags and one black. We have $l(\tau )=0$
if and only if the vertex $v_0$ carrying $a$ is not critical.
We have $l(\tau )=n\ge 1$ if and only if
$v_0$ is critical, and there is a (unique) sequence
of pairwise distinct vertices $v_0,\dots ,v_n:=w$ such that
$v_1,\dots ,v_{n-1}$ are critical, whereas $v_n$ is not,
and that $v_i,v_{i+1}$ are neighbors, that is ends of an edge.

\smallskip

We will sometimes call $l(\tau )$ simply length of $\tau$.

\smallskip

We will say that $w$ (and $\tau$) is of type I, if
$w$ carries only two white flags, and therefore
$\ge 2$ black flags. We will say that $w$ (and $\tau$) 
is of type II, if
$w$ carries $\ge 3$ white flags.

\medskip

The $n$--th step of induction, $n\ge 0$, will consist of
the following constructions and verifications.

\smallskip

(A)${}_n$. {\it A definition of $\rho_*^H (m(\tau ))$ for all 
$S$--labeled stable trees $\tau$
with $l(\tau )=n.$}

\smallskip

Actually, this definition generally will depend on arbitrary choices
and produce directly only an element of $\oplus k m(\sigma )$ where $\sigma$
runs over $S'$--labeled painted stable trees.
Hence we will have to check that

\smallskip

(B)${}_n$. {\it For $l(\tau )= n,$ $\rho_*^H (m(\tau ))$ 
are defined unambiguously
modulo $I_{S'}.$}

\smallskip

Finally, we will have to check that $\rho_*^H (m(\tau ))$
depends only on $[m(\tau )]\in H^*_S$, or equivalently

\smallskip

(C)${}_n$. {\it $\rho_*^H$ extended by linearity sends to $I_{S'}$
each standard relation in $I_S$ whose all terms have length $\le n$.}

\medskip

{\bf 3.9. Construction of $\rho^H_*$: the case $l(\tau )=0$.}

(A)${}_0$. If $\tau$ is of type I, we put $\rho_*^H(m(\tau ))=0$.
The same prescription will hold for type I and any length.

\smallskip

If $\tau$ is of type II, repainting the label $a$ produces
a stable $S'$--tree, say, $\tau'$. In this case we put
$\rho_*^H(m(\tau ))=[m(\tau')].$

\smallskip

(B)${}_0$. Clearly, this prescription is unambiguous.

\smallskip

(C)${}_0$. Let $R(\sigma ,u; I,J,K,L)$ be a relation (3.8),
such that all its terms $\sigma (\alpha )$ are of length 0.
Since $\sigma$ is obtained from any $\sigma (\alpha )$
by collapsing an edge, we have $l(\sigma )=0.$ If $u\ne w=v_0$,
the type of $\sigma$ is the same as the type of all
$\sigma (\alpha )$, and we get either
$\rho(R(\sigma ,u; I,J,K,L))=0$ (type I), or
$\rho(R(\sigma ,u; I,J,K,L))=R(\sigma' ,u; I,J,K,L)$
(type II). If $u=w=v_0$ and $\sigma$ is of type I,
all $\sigma (\alpha )$ must be of type I, and  again
$\rho(R(\sigma ,u; I,J,K,L))=0$.

\smallskip

It remains to consider the case when $u=w=v_0$
and $\sigma$ is of type II. If $a\notin \{I,J,K,L\}$,
then for any stable partition $\alpha$
of $F_{\sigma}(u)$, $a$ gets into the part of $\alpha$
containing an extra white label from $\{I,J,K,L\}$.
Hence $l(\sigma (\alpha ))=0$ and moreover,
$\sigma (\alpha )$ is of type II. Therefore
$\rho(R(\sigma ,u; I,J,K,L))=R(\sigma' ,u; I,J,K,L))$.
Finally, let $a\in \{I,J,K,L\}$, say, $a=I$
so that $I$ is white. We assumed that all terms
of $R(\sigma ,u; I,J,K,L)$ are of length $0$.
Hence $J$ and $L$ must be white as well,
and we have again 
$\rho(R(\sigma ,u; I,J,K,L))=R(\sigma' ,u; I,J,K,L))$.

\medskip

{\bf 3.9.1. Construction of $\rho_H^*$.} Let
$\tau'$ be a painted stable $S'$--tree.
The label $a$ in it is black, so let us
denote by $\tau$ the result of repainting it white.
Then $l(\tau )=0.$ We put
$\rho_H^*(m(\tau')):=[m(\tau )].$ Clearly,
after such reverse repainting, each standard
relation becomes a standard relation,
so that we have a well defined map
$\rho_H^*:\,H^*_{S'}\to H^*_S.$ The commutativity
of (3.12) is obvious.
Since this
construction involves only $S$--trees of length $0$,
the discussion above already shows that 
$\rho^H_*\circ \rho_H^* =\roman{id}$ identically.

\smallskip

We will now return to the direct image. From now on,
we assume that $n\ge 1$ and that the statements
(A)${}_m$, (B)${}_m$, (C)${}_m$ have been already treated
for all $m\le n-1.$

\medskip

{\bf 3.10. Passage from $n-1$ to $n$: prescription (A)${}_n$.}
Let $l(\tau )=n$, and $\tau$ be of type II
(for type I, repainting produces zero for any length). 
Denote by $\sigma$ the result of collapsing to a vertex $u$
the last edge $e_{n-1}$ with vertices
$v_{n-1},v_n$ of the critical branch of $\tau$. Since $w=v_n$ in $\tau$
carried $\le 3$ white flags, $u$ carries
$\le 2$ white flags, say, $K,L$, besides
the white flag, say, $I$, which starts the way from $u$ to
$a$ for $n\ge 2$, or coincides with $a$ for $n=1$.
Denote by $J\in F_u(\sigma )$ the unique black flag
that was carried by the critical vertex $v_{n-1}.$

\smallskip

With this notation, consider the relation $R(\sigma, u; I,J,K,L).$
It contains exactly one term $\sigma (\alpha_0)$
of length $n$, namely for $\alpha_0=IJ|KL...$.
We have $\sigma(\alpha_0)=\tau$, and $l(\sigma(\alpha ))=n-1$
for $\alpha \ne \alpha_0.$ By the inductive assumption, (the lifts of)
$\rho^*_H (\sigma(\alpha ))$ are well defined
modulo the subspace generated by the standard relations
all terms of which have length $\le n-1.$
We put
$$
\rho_*^H (\tau ):= -\sum_{\alpha\ne \alpha_0:\,IJ\alpha KL}
\rho_*^H (\sigma(\alpha )) +
\sum_{\alpha :\,KJ\alpha IL}
\rho_*^H (\sigma(\alpha ))\, \roman{mod}\, I_{S'}.
\eqno(3.13)
$$
Equivalently, repainting
of such $R(\sigma, u; I,J,K,L)$ lands in $I_{S'}.$ Thus, applying
the pescription (A)${}_n$ is the same as postulating
some paticular cases of the statement (C)${}_n$. 

\medskip

{\bf 3.11. Passage from $n-1$ to $n$: the statement (B)${}_n$.} 
We have to check that (3.13) does not depend on the
arbitrary choices of $K,L$. To pass from one couple 
of white flags $K,L$ to another
one it suffices to replace one flag in turn, so we will consider two cases.

\smallskip

{\it Replacement $K\mapsto K'$.} We know that after repainting
$R(\sigma, u; I,J,K,L)$ we land in $I_{S'}$ and wish to establish
the same for $R(\sigma, u; I,J,K',L)$. A version of (3.10) gives
$$
R(\sigma, u; I,J,K,L) - R(\sigma, u; I,J,K',L) =
R(\sigma, u; K,J,K',L)\, .
\eqno(3.14)
$$
Two terms of length $n$ cancel in the left hand side of (3.14).
Hence $R(\sigma, u; K,J,K',L)$ contains only terms
of length $\le n-1$, and the result follows by application
of (C)${}_{n-1}$.

\medskip

{\it Replacement $L\mapsto L'$.} A similar
calculation shows that the respective difference will be now 
$$
R(\sigma , u; L',K,L,I)
\eqno(3.15)
$$
again with all terms of length $n-1$ and the same conclusion. 

\smallskip

For further use, notice that in this reasoning
$L'$ might have been black as well: all involved quadruples
remain allowed. Therefore
any relation of the type $R(\sigma , u; I,J,K,L)$
with white $I,K$ and black $J,L$ will be also
repainted to an element of $I_{S'}.$

\smallskip

It remains to prove (C)${}_n$. Before doing this,
it is convenient to review the structure
of the standard relations that must be treated
at the $n$--th step.

\medskip

{\bf 3.12. Standard relations having terms of length $n$.}
Let $R(\sigma ,u;I,J,K,L)$ be a relation
all terms of which have length $\le n$, and such that this
bound is achieved.

\smallskip

Consider a term $\sigma (\alpha_0)$ of length $n$.
Let $v_0,\dots ,v_{n-1},v_n=w$ be (the sequence of vertices of) the
critical branch of $\sigma (\alpha_0)$ so that 
$a$ is incident to $v_0$, all vertices
$v_0,\dots ,v_{n-1}$ are critical, and
$w$ is not critical. Let $e_i$ be the edge
with vertices $v_i,v_{i+1}.$

\smallskip

Consider the position of the edge $e(\alpha_0)$ in $\sigma (\alpha_0)$
which was created by $\alpha_0$ and gets collapsed into $u$
in $\sigma$. We will have to treat separately
the following (exhaustive) list of alternatives (a), (b), (c), (d).

\smallskip

(a) {\it $e(\alpha_0 )=e_i$ for some $i\le n-2$.}

\smallskip

This is, of course, possible only for $n\ge 2$.

\smallskip

In this case $F_{\sigma}(u)$ consists of four flags:
white flag leading from $u$ to $a$ (or $a$ itself
for $i=0$) which we denote $I$; black flag
$J$ incident to $v_i$ in $\sigma (\alpha )$, black flag $L$ incident to $v_{i+1}$ in $\sigma (\alpha )$,
white flag $K$ leading from $u$ to $w$.
With this notation, $\alpha_0 = IJ|KL.$ Put
$\overline{\alpha}_0:= IL|JK.$
Up to renaming flags and changing signs,
we get a relation comprising only two nonvanishing
terms, both of length $n$.
$$
R(\sigma ,u; I,J,K,L)=m(\sigma (\alpha_0))-m(\sigma (\overline{\alpha}_0 )).
\eqno(3.16)
$$
We can call (3.16) ``exchange of two neighboring
black tails on the critical branch of $\sigma .$''
Notice that the length of $\sigma$ is $i-1$.

\smallskip

(b) {\it $e(\alpha_0 )=e_{n-1}$ so that $w$ is a vertex of
$e(\alpha_0).$}

\smallskip

We may and will define $I,J$ as above. If $\sigma (\alpha_0)$ is of type
I,  then $\sigma$ and all $\sigma (\alpha )$ are of type I 
so that
$R(\sigma ,u;I,J,K,L )$ must be repainted to zero.
Type II will be treated below.  

\smallskip

(c) {\it None of the above, but 
$w$ is a vertex of $e(\alpha_0 ).$} 

\smallskip

This will be the most difficult case.

\smallskip

(d) {\it None of the vertices $v_0,\dots , v_n=w$
is a vertex of $e(\alpha_0 ).$}

\smallskip

In this case all terms $\sigma (\alpha )$ as well
as $\sigma$ have length $n$ and in fact
share the common critical branch. 

\smallskip

We will now treat these options in the following
order: (b), (a), (d), (c).

\medskip 

{\bf 3.13. Passage from $n-1$ to $n$: the statement (C)${}_n$ in the case (b).}
We use the notation defined in 3.12 (b) and assume
that $\sigma$ is of type II. Since
$\{I,J,K,L\}$ is allowable, $I$ is white, and
$J$ is black, $K$ must be white. 

\smallskip

If $L$ is white as well, then $\sigma (\alpha_0)$
is the only term of length $n$ in $R(\sigma ,u;I,J,K,L)$,
and the whole relation repaints to an element of $I_{S'}$:
cf. the last remark in 3.10.

\smallskip

If $L$ is black, we refer to the remark
made after formula (3.15) which proves the same result. 

\medskip 

{\bf 3.14. Passage from $n-1$ to $n$: the statement (C)${}_n$ in the case (a).}
Again, we may assume that $\sigma (\alpha_0 )$ is
of type II.
In order to repaint (3.16), we have to apply
to both terms the prescription spelled out
in 3.10, formula (3.13). However, the notation
adopted in (3.13) conflicts with the one
adopted in (3.16). We will keep (3.16)
and rewrite (3.13) as follows. Let $I'$
be the white flag at $w$ in $\sigma (\alpha_0)$
and in $\sigma (\overline{\alpha}_0)$
leading to $a$, $J',K',L'$ three other
flags at $w$ such that $K',L'$ are white.
Then we can use $R(\sigma (\alpha_0),u;I',J',K',L')$
and the similar relation for $\overline{\alpha}_0$
in order to replace $m(\sigma (\alpha_0))$
and $m(\sigma (\overline{\alpha}_0 ))$ modulo $I_S$
by  linear
combinations of good monomials corresponding to trees
of length $\le n-1$.

\smallskip

If $i+1< n-1$, we can then pair the terms in the resulting
difference in such a way that we will get
a linear combination of the exchange relations of the type (3.16)
written for a new family of trees which are of smaller length.

\medskip

{\bf 3.15. Passage from $n-1$ to $n$: the statement (C)${}_n$ in the case (d).}
In this case there is no
interaction between the surgeries made upon $\sigma$
by the partitions $\alpha$ involved
in $R(\sigma , u;I,J,K,L)$ and the partitions
$\beta$ of $F_{\sigma}(w)$ which occur in a 
formula of the type (3.13) which can be chosen common
for all terms $\sigma (\alpha )$. Hence we can
reorder the surgeries and the summations and start with summing
over $\beta$. This will show
that $R(\sigma , u;I,J,K,L)$ is repainted to a sum
of standard relations. 

\medskip

{\bf 3.16. Passage from $n-1$ to $n$: the statement (C)${}_n$ in the case (c).} Consider again
a term $\sigma (\alpha_0)$ of length $n$. The edge created by $\alpha_0$
in  $\sigma (\alpha_0)$ is incident to the last vertex of
the critical branch, and another vertex of this edge does not
belong to the critical branch. This edge collapses to the vertex $u$ of
$\sigma$. Denote by $M$ the flag at $u=w=v_n$ (notation now refers to $\sigma$)
belonging to the critical branch. Thus $M$ leads to $v_{n-1}$. Denote by
$B$ the (single) black flag at $v_{n-1}$ and by $A$
the single white flag at $v_{n-1}$ leading in the direction of $a$ (if $n=1$,
$B:=b$, $A:=a$.) Furthermore, we have an allowed
quadruple of flags $\{I,J,K,L\}\in F_{\sigma}(u)$
defining the relation $R(\sigma ,u;I,J,K,L)$ that we are going to reoaint.

\medskip

{\bf 3.16.1. Claim.} {\it The total statement (C)${}_n$ in the case (c)
will follow if we prove it assuming the following condition:

\smallskip

(*) The critical branch of $\sigma$ is of length $n$,  $I=M$, $J,L$ are white.}

\smallskip

{\bf Proof of the Claim.} In the notation explained at the beginning of 3.16,
we have either $M\notin \{I,J,K,L\},$ or $M\in \{I,J,K,L\}.$

\smallskip

If $M\notin \{I,J,K,L\},$ we may and will assume that $I,K$ are white.
In fact, from (3.9) it follows that this can be achieved by renaming 
the flags and
changing the sign of the relation if need be. Since for any $\sigma (\alpha )$
in $R(\sigma ,u;I,J,K,L)$ the flag $M$ gets at the same vertex as
either $I$, or $K$, its vertex remains non--critical.
Therefore $l(\sigma (\alpha ))=l(\sigma )=n$, and all $\sigma (\alpha )$
share the common critical branch with $\sigma$.

\smallskip

Now write the relation (3.10):
$$
R(\sigma ,u;I,J,K,L)=R(\sigma ,u;M,K,J,I)+R(\sigma ,u;M,I,L,K)
$$
In both groups of flags $M,K,J,I$ and $M,I,L,K$ appearing at the right
hand side   the first, second, and fourth terms
are white. Since $I,K$ are white, the vertex of $M$ in each term 
of the relations at the right hand side is not critical. This means
that both these relations satisfy the condition (*), and if (C)${}_n$
holds for them, it holds for $R(\sigma ,u;I,J,K,L)$ as well.

\smallskip

Let now $M\in \{I,J,K,L\}.$ We may and will assume that $M=I$ 
so that $I$ is white. If moreover both $J$ and $L$ are white, 
then (*) holds as above.

\smallskip

Let us show that cases when $J$ or $L$ is black fall into 
another category. In fact, if say $J$ is black,
then the partition $\alpha_0:=IJ|KL ...$ creates 
a new critical edge in $\sigma (\alpha_0)$ so that we are in the case
3.12(a) or 3.12(b) treated earlier.

\medskip

{\bf 3.16.2. Treatment of the case 3.16.1 (*).} Put $F:=F_{\sigma}(u) -\{I,J,K,L\}.$
To see better how $R(\sigma ,u;I,J,K,L)$ repaints to zero modulo $I_{S'}$
consider first the simplest case in which $F=\emptyset$. 
We keep notation explained at the beginning of 3.16. The last edge
of the critical branch of $\sigma$ has two vertices
carrying flags $A,B$, resp. $J,K,L$, besides halves of the edge itself.
We will denote $\sigma$ symbolically by $AB|JKL$. There will be two
partitions determining terms $\sigma (\alpha )$ of $R(\sigma ,u;I,J,K,L)$;
the respective trees will be denoted $AB|J|KL$ and $AB|L|KJ$ so that
in the current shorthand notation
$$
R(\sigma ,u;I,J,K,L) = AB|J|KL - AB|L|KJ \, .
\eqno(3.17)
$$
Applying the prescription (3.13) for repainting two terms
of (3.17),  we have a choice: either $A$, or $B$ can be moved
to the middle vertex. We decide to move $A$,  interchanging $A$ with $J$, resp. with $L$.
At this step the assumption 3.16.1 (*) that $J$ and $L$ are white
is critically used: otherwise one or both trees in (3.17)
would become unstable. We  get symbolically
$$
\rho_*^H (R(\sigma ,u;I,J,K,L)) = \rho_*^H (BJ|A|KL) - \rho_*^H (BL|A|JK) \, \roman{mod}\,
I_{S'}\,  
\eqno(3.18)
$$ 
where at the right hand side both trees have  length $n-1$.
Now we will in three consecutive steps interchange $A$ with $K$, $B$ with $K$,
$B$ with $A$. More formally, add to (3.18) appropriate
relations of length $n-1$ which are repainted into $I_{S'}$ in view of
(C)${}_{n-1}$. For example, interchanging $A$ with $K$
in the first term of the right hand of (3.18) means subtracting
(the repainting of) the standard relation $BJ|A|KL - BJ|K|AL$.

\smallskip

Thus, omitting $\rho_*^H$ for brevity, we rewrite (3.18) modulo $I_{S'}$ 
consecutively into 
$$
BJ|K|AL - BL|K|AJ,
\eqno(3.19)
$$
then
$$
KJ|B|AL - KL|B|AJ,
\eqno(3.20)
$$
and finally
$$
KJ|A|BL - KL|A|BJ.
\eqno(3.21)
$$
Now a miracle happens: $KJ|A|BL$ as a tree is isomorphic to $BL|A|JK$,
and $KL|A|BJ$ to $BJ|A|KL$. Therefore (3.21) differs from
(3.18) by the sign, and hence both expressions vanish
over any ring where $2$ is invertible.

\smallskip

We now take a deep breath preparing for the last stretch of the proof,
and consider the case of arbitrary $F:=F_{\sigma}(u) -\{I,J,K,L\}.$

\smallskip

The formula (3.17) must now be replaced by the sum taken over all
{\it ordered} 2--partitions  $F=F(0)\cup F(1)$ (corresponding to former $\alpha$'s):
$$
R(\sigma ,u;I,J,K,L) = \sum_{F(0),F(1)} \left[\, AB|JF(0)|F(1)KL - AB|LF(0)|F(1)JK 
\,\right]\, 
\eqno(3.22)
$$
Here we write say, $F(1)KL$, in place of $F(1)\cup \{K,L\}$. Moreover, we interpret
say, the expression $AB|JF(0)|F(1)KL$ as a notation for the
isomorphism class of the following $S$--tree: take a linear tree with two edges and
three consecutive vertices $u_1,u_2,u_3$  and attach to the vertices the following
branches of the tree $\sigma$ which are denoted by their initial flags:
$A$ and $B$ to $u_1$, $\{J\}\cup F(0)$ to $u_2$, $F(1)\cup \{KL\}$ to $u_3$.

\smallskip

After repainting $a$ black, each term of (3.22) turns into a similar sum,
now taken over partitions $F(0)=F(00)\cup F(01)$,
(we omit $\rho^*_H$ for brevity):
$$
AB|JF(0)|F(1)KL \cong 
$$
$$
-\sum \Sb F(00), F(01)\\ F(00)\ne \emptyset \endSb
  ABF(00)|JF(01)|F(1)KL 
+\sum \Sb F(00), F(01) \endSb
  BJF(00)|AF(01)|F(1)KL \, ,
\eqno(3.23a)
$$
$$
-AB|LF(0)|F(1)JK \cong 
$$
$$
\sum \Sb F(00), F(01)\\ F(00)\ne \emptyset \endSb
  ABF(00)|LF(01)|F(1)JK
-\sum \Sb F(00), F(01) \endSb
  BLF(00)|AF(01)|F(1)JK \, .
\eqno(3.23b)
$$
The first sums at the right hand sides of (2.23a), (2.23b)
did not appear in  the degenerate case (3.18) where $F$ was empty, 
so they must be treated separately. Let us group
the respective terms together in the following way. 
$$
\sum \Sb F(00)\subset F\\F(00)\ne \emptyset \endSb
\sum \Sb F(01), F(1)\\F(01)\cup F(1) =F-F(00) \endSb
\left[\,ABF(00)|LF(01)|F(1)JK-  ABF(00)|JF(01)|F(1)KL\,\right] \, .
$$
Each inner sum here is (the result of repainting of) a
standard relation. Each term in all these relations has length $\le n-1$
because $F(00)\ne \emptyset$. Hence all these terms can be disposed off
thanks to (C)${}_{n-1}.$

\smallskip

Turning now to the last sums in (3.23a), (3.23b),
we consecutively transform them on the pattern of (3.19)--(3.21).

\smallskip

We describe the first step corresponding to the interchange
of $A$ and $K$ in some detail. Sum second terms in (3.23a), resp. (3.23b)
over all $F(0)\subset F$ and then make the summation over $F(00)\subset F$
external. Each inner sum will become ``one half''
of a standard relation, and we replace it by another half.
To write a formula, extending (3.19), we denote in the inner sum, 
when $F(00)$ is fixed,
$G:=F-F(00)=F(01)\cup F(1)$. The inner sum
will be extended over partitions $G=G(0)\cup G(1)$. The general case
of (3.19) takes form
$$
\sum_{F(00)} \sum_{G(0),G(1)} \left[\,BJF(00)|KG(0)|ALG(1) -
BLF(00)|KG(0)|AJG(1)\,\right]\, .
\eqno(3.24)
$$
Now reshuffle again, make summation over $G(1)$ external,
and with fixed $G(1)$ put $H=F-G(1)=F(00)\cup G(0).$
The general case
of (3.20) takes form
$$
\sum_{G(1) } \sum_{H(0),H(1)} \left[\,KJH(0)|BH(1)|ALG(1) -
KLH(0)|BH(1)|AJG(1)\,\right]\, .
\eqno(3.25)
$$
Finally, at the last step make summation over $H(0)$ external and put $E=F-H(0)$
so that (3.21) becomes
$$
\sum_{H(0)} \sum_{E(0),E(1)} \left[\,KJH(0)|AE(0)|BLE(1) -
KLH(0)|AE(0)|BJE(1)\,\right]\, .
\eqno(3.26)
$$
As above, (3.26) differs by sign from the sum
of all last terms in (3.23a) and (3.23b), because
the triples $(H(0), E(0), E(1))$ in (3.26) and $(F(1),
F(01),F(00))$ in (3.23) run over all
ordered partitions of $F$ into three pairwise disjoint subsets.

\smallskip

This completes the construction of $\rho_*^H$ and the proof
of the Theorem 3.1.1.

\bigskip

\centerline{\bf References}

\medskip

[H] B.~Hassett. {\it Moduli spaces of weighted pointed stable curves.}
Preprint \newline math.AG/0205009

\smallskip

[Ke] S.~Keel. {\it Intersection theory of moduli space
of stable $N$--pointed curves of genus zero}.
Trans.~AMS, 330:2 (1992), 545--574.

\smallskip

[Kn] Knudsen F.F. {\it The projectivity of the moduli space
    of stable curves II.
The stacks $\overline M_{0,n}$}. Math. Scand. 52 (1983), 163--199.

\smallskip

[KoMa] M.~Kontsevich, Yu.~Manin. {\it Gromov--Witten classes, quantum cohomology, and enumerative 
geometry.} Comm. Math. Phys.,
164:3 (1994), 525--562.

\smallskip 

[KoMaK] M.~Kontsevich, Yu.~Manin (with Appendix by
R.~Kaufmann). {\it Quantum cohomology of a product.} Inv. Math., 124, f. 1--3 (Remmert's Festschrift) (1996), 313--339.

\smallskip

[LoMa1] A.~Losev, Yu.~Manin. {\it New moduli spaces of pointed 
curves and pencils of flat connections.}
Michigan Journ. of Math., vol. 48
(Fulton's Festschrift), 2000, 443--472. Preprint math.AG/0001003

\smallskip

[LoMa2] A.~Losev, Yu.~Manin. {\it Extended modular operad.} 
Preprint math.AG/0301003

\smallskip

[Ma] Yu.~Manin. {\it Frobenius manifolds, quantum cohomology,
and moduli spaces.}  AMS Colloquium Publications, vol. 47, Providence, RI, 1999, xiii+303 pp.

\enddocument